\documentclass[11pt]{article}
\usepackage[utf8]{inputenc}
\usepackage{amsmath, amsfonts}
\usepackage{tikz}

\newcommand{\qtt}{Q_{2\tau}}
\newcommand{\norm}[1]{\left|\left| #1 \right|\right|}

\title{An Evaluation of novel method of Ill-Posed Problem for the Black-Scholes Equation solution }
\author{General Paper}
\author{Kirill V.Golubnichiy, Tianyang Wang and Andrey V. Nikitin
\and Department of Mathematics, 
\and Paul G. Allen School of Computer Science \& Engineering
\and University of Washington, Seattle, USA
\and kgolubni@math.washington.edu, tianyang@uw.edu  
\and andrey.nikitin@hotmail.com}

\date{}

\usepackage[square]{natbib}
\usepackage{graphicx}
\counterwithout{figure}{section}
\counterwithout{figure}{subsection}

\setcitestyle{numbers}

\begin{document}

\maketitle

 \textbf{\text{Abstract}}

It was proposed in \cite{KlibGol} a new empirical mathematical method to work with the Black-Scholes equation. This equation is solved  forwards in time to forecast prices of stock options. It was used the regularization method because of ill-posed problems. Uniqueness, stability and convergence  theorems for this method are formulated. For each individual option, historical data is used for input. The latter is done for two hundred thousand stock options selected from the Bloomberg terminal \cite{Bloom} of University of Washington. It used the index Russell 2000 \cite{Russ}. The main observation is that it was demonstrated that technique, combined with a new trading strategy, results in a significant profit on those options. On the other hand, it was demonstrated the trivial extrapolation techniques results in much lesser profit on those options. This was an experimental work. The minimization process was performed by Hyak Next Generation Supercomputer of the research computing club of University of Washington. As a result, it obtained about 50,000 minimizers \cite{Git}. The code is parallelized in order to maximize the performance on supercomputer clusters. Python with the SciPy module was used for implementation. You may find minimizers in the source package that is available on Github \cite{Git}. Chapter 7 is dedicated to application of machine learning. We were able to improve our results of profitability using minimizers as new data. We classified the minimizer's set to filter for the trading strategy.  All results are available on \cite{Git}.   

\vspace{0.5em}
 \textbf{\text{Notations}}

\begin{itemize}
    \item $\qtt \equiv (s_b, s_a) \times (0, 2\tau) \subset \mathbb{R}^2$.
    \item $\partial \qtt \equiv \{(s, t): t = 0, (s = s_a \text{ or } s = s_b)\}$.
    \item $\mathcal{L}: H^{2, 1}(\qtt) \rightarrow L^2(\qtt)$ is the differential operator, where
    \[
    \mathcal{L} u \equiv \frac{\partial u}{\partial t} + \frac{\sigma^2(t)}{2} s^2 \frac{\partial^2 u}{\partial s^2}.
    \]
    \item Hom$(A, B)$ is the set of all homomorphic functions $f: A \rightarrow B$. We didn't use the convention $L(A, B)$ in order to distinguish from the matrix $L = D_t + RD_{ss}$.
    \item $L^2(\qtt)$ space is the space of measurable functions for which the 2-nd power of the absolute value is Lebesgue integrable, and the norm is defined as:
    \[\norm{f}_{L^2(\qtt)} \equiv \left( \int_{\qtt} \left| f \right|^2 d\mu \right)^{1/2} < \infty. \]
    \item $J_\beta(u) \equiv \int_{\qtt} (\mathcal{L}u)^2 dsdt + \beta \norm{u - F}_{L^2(\qtt)}$ is the Tikhonov-like functional. 
    \item $\hat{J}_\beta(u)$ is the discrete version of $J_\beta(u)$ (derived by T-shape finite difference scheme).
    \item $F(s, t).$
\end{itemize}
 
  \textbf{\text{Keywords:}}
  
  Black-Scholes equation, Ill-posed problem, regularization method, parabolic equation with the reversed time,  T-shape, machine learning, neural network.  

\section{Introduction}
In 1943 Tikhonov invented the regularization method for ill-posed problems \cite{Tikhonov1943Stability}. He introduced the fundamental concept. This concept consists of the following
three conditions which should be in place when solving the ill-posed problem.

\bigskip

1. One should a priori assume that there exists an ideal exact solution $x*$ of the problem for an ideal noiseless data $y$.

2. The correctness set should be chosen a priori, meaning that some a priori
bounds imposed on the solution $x$ of the problem should be imposed.

3. To construct a stable numerical method for the problem.

\bigskip

For the first time, the question about global uniqueness theorems was addressed
positively and for a broad class of Coefficient Inverse Problems with single measurement data in the
works of A.L. Bukhgeim and M.V. Klibanov in 1981. These results were published in their paper \cite{KlibanovBuk}. After widely applied to many physical examples. The first complete proofs were published in
two separate papers \cite{Buk} and \cite{Klib1}. This technique
is now called the “Bukhgeim-Klibanov method.” This method is the
only one enabling for proofs of global uniqueness results for multidimensional
Coefficient Inverse Problems with single measurement data. The Bukhgeim-Klibanov method is based on the idea of applications of the so-
called Carleman estimates to proofs of uniqueness results for Coefficient Inverse Problems. The main interest in applications in, for example, the hyperbolic case, is when one of
initial conditions is identically zero and another one is either the $\delta$ function or that
the wave field is initialized by the plane wave. The uniqueness question in the latter
case remains a long-standing and well-known unsolved problem.

\bigskip

\vspace{1em}
\section{The mathematical model}

Find an approximate solution of the Black-Scholes equation%
\begin{equation}
Lu=u_{t}+\frac{\sigma ^{2}\left( t\right) }{2}s^{2}u_{ss}=0,\left(
s,t\right) \in \left( s_{b}\left( 0\right) ,s_{a}\left( 0\right) \right)
\times \left( 0,2\tau \right) =Q_{2\tau },
\end{equation}%
subject to boundary conditions \begin{equation}
u\left( s_{b},t\right) =u_{b}\left( t\right) ,u\left( s_{a},t\right)
=u_{a}\left( t\right) ,\text{ }t\in \left[ 0,2\tau \right] ,
\end{equation}
and the initial condition
\begin{equation}
u\left( s,0\right) =f\left( s\right) ,\text{ \ \ }s\in \left[ s_{b}\left(
0\right) ,s_{a}\left( 0\right) \right] .
\end{equation}%

$s$ is the stock price, $t$ is time, $\sigma \left( t\right) $ is the
volatility of the stock option 

$u\left( s,t\right) $ is the price of the stock option.

We predict option price from ``today" to ``tomorrow" and ``the day after tomorrow". 255 trading days annually. ``One day" $\tau =1/255.$ ``Today" $t=0.$ ``Tomorrow" $t=\tau .$ ``The day after tomorrow" $t=2\tau .$ $s-$interval: $s\in \left[ s_{b}\left( 0\right) ,s_{a}\left( 0\right) \right].$

\bigskip

To solve the problem, we minimize
the functional $J_{\beta }\left( u\right) $,%
\begin{equation*}
J_{\beta }\left( u\right) =\int_{Q_{2\tau }}\left( Lu\right)
^{2}dsdt+\beta \left\Vert (u-F)\right\Vert _{H^{2}\left( Q_{2\tau }\right) }^{2},
\end{equation*}
where $\beta \in \left( 0,1\right) $ is the regularization parameter and 
$F(s,t)=s(u_{a}(t)-u_{b}(t))+u_{b}(t),$ $(s,t)\in Q_{2\tau}.$  The function $F \in H^{2}(Q_{2\tau}).$ It follows from (2) and (3) that 

\begin{equation}
F\left( s,0\right) =f\left( s\right) ,  \label{4000}
\end{equation}%

\begin{equation}
F\left( s_{b},t\right) =u_{b}\left( t\right) ,F\left( s_{a},t\right)
=u_{a}\left( t\right).  \label{3000}
\end{equation}

Uniqueness and existence of the minimizer $u_{\beta }\in H^{2}\left(
Q_{2\tau }\right) $ follow from the Riesz theorem.

Convergence of minimizers to the exact solution when the level of error in
the boundary and initial data tends to zero, was proven by Klibanov (2015) \cite{KlibGol}
using a Carleman estimate.

\vspace{1em}
\section{T-shape finite difference scheme}

\subsection{Discretization}

The domain of the PDE is given by $\qtt \equiv (s_b, s_a) \times (0, 2\tau) \subset \mathbb{R}^2$. The domain can be discretizied by tuples $(s, t)$: \cite{Wolf} $s_i = s_b + i \cdot ds$, $t_j = j \cdot dt$, for $i, j = 0, ..., M$. The actual value of $M$ can be determined depending on the performance of the computer and the required accuracy; $M \geq 3$ must be true otherwise it is impossible to solve, and we noticed that the results converges well (relative difference between different $M$'s $< 0.01$) when $M$ is greater than 20.

Let $u_{(i, j)}$ denote the option price corresponding to a stock price $s_i$ at time $t_j$: $u_{(i, j)} = u(s_i, t_j)$ \cite{Wolf}. 

For discretizing the partial derivatives, we considered $\frac{\partial^2 u}{\partial s^2}$ and $\frac{\partial u}{\partial t}$ separately. We used the backward difference scheme for the first derivative, and the standard central difference scheme \cite{Strikwerda} for the second derivative, i.e.

\begin{equation}
    \delta_{t-} u_{(i, j)} = \frac{u_{(i, j)} - u_{(i, j-1)}}{dt};
\end{equation}
\begin{equation}
    \delta^2_{s} u_{(i, j)} = \frac{u_{(i+1, j)} - 2u_{(i, j)} + u_{(i-1, j)}}{ds^2}.
\end{equation}

Then the partial differential equation can be approximated by:
\begin{equation}
        \delta_{t-} u_{(i, j)} + \frac{\sigma^2(t_j)}{2} s_i^2 \delta^2_{s} u_{(i, j)} = 0.
\end{equation}

Consider the equation above at a single point $u_{(i,j)}$. It involves four adjacent points, including $u_{(i, j)}$ itself: $u_{(i, j)}, u_{(i-1, j)}, u_{(i+1, j)}, u_{(i, j-1)}$ . The four points form a letter ``T'' shape in the Cartesian coordinates:

\begin{figure}
\begin{center}
\begin{tikzpicture}
\node (0) at (0.8, -0.6) {$(i, j)$};
\node (1) at (3, -0.6) {$(i+1, j)$};
\node (1) at (-1, -0.6) {$(i-1, j)$};
\node (1) at (1, -2.6) {$(i, j-1)$};
\filldraw (0, 0) circle (0.1);
\draw (0, 2) circle (0.1);
\draw (0, -4) circle (0.1);
\filldraw (0, -2) circle (0.1);
\filldraw (2, 0) circle (0.1);
\draw (2, 2) circle (0.1);
\draw (2, -4) circle (0.1);
\draw (2, -2) circle (0.1);
\draw (4, 0) circle (0.1);
\draw (4, 2) circle (0.1);
\draw (4, -4) circle (0.1);
\draw (4, -2) circle (0.1);
\filldraw (-2, 0) circle (0.1);
\draw (-2, 2) circle (0.1);
\draw (-2, -4) circle (0.1);
\draw (-2, -2) circle (0.1);
\draw (-4, 0) circle (0.1);
\draw (-4, 2) circle (0.1);
\draw (-4, -4) circle (0.1);
\draw (-4, -2) circle (0.1);

\draw (-4, 2) -- (4, 2);
\draw (-4, 0) -- (4, 0);
\draw (-4, -2) -- (4, -2);
\draw (-4, -4) -- (4, -4);
\draw (-4, 2) -- (-4, -4);
\draw (-2, 2) -- (-2, -4);
\draw (0, 2) -- (0, -4);
\draw (2, 2) -- (2, -4);
\draw (4, 2) -- (4, -4);
\draw[very thick] (-2, 0) -- (2, 0);
\draw[very thick] (0, 0) -- (0, -2);
\end{tikzpicture}
\end{center}
\caption{The four grid points involved in the T-shape}
\end{figure}
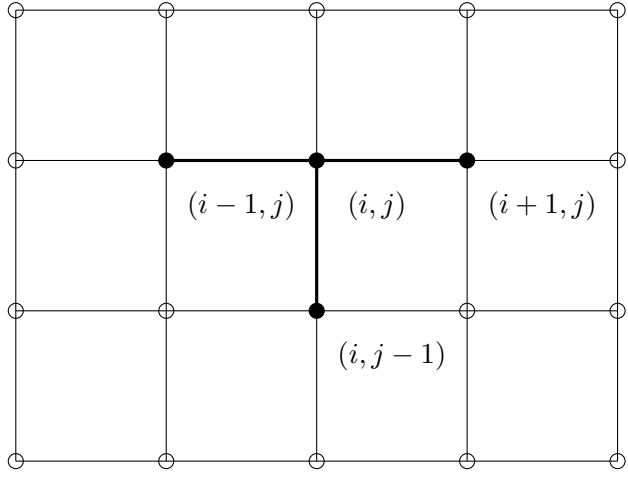

Therefore we called it \textbf{T-shape finite difference scheme}. Denote $\partial \hat{Q}_{2\tau} = \{(i, j): u_{(i,j)} \in \partial\qtt\}$ and $\hat{Q}_{2\tau} = \{(i, j): u_{(i,j)} \in \qtt\}$; as a result, the Tikhonov-like functional $J_\beta (u)$ \cite{KlibGol} introduced above can also be approximated by:

\begin{multline}
    \hat{J}_\beta(u) = dtds \sum_{(i, j) \in \hat{Q}_{2\tau}\setminus\partial \hat{Q}_{2\tau}} \left( \delta_{t-} u_{(i, j)} + \frac{\sigma^2(t_j)}{2} s_i^2 \delta^2_{s} u_{(i, j)} \right)^2  \\
     + \beta dtds \sum_{(i, j) \in \hat{Q}_{2\tau}} \left( u_{(i, j)} - F(s_i, t_j)\right)^2.
\end{multline}

\subsection{Serialization}

In order to make the minimization of $\hat{J}_\beta(u)$ ``code friendly'', we need to find a better way to represent the functional instead of a big summation. Notice that the T-shape finite difference scheme, when regarding the four points as a vector $v \in \mathbb{R}^4$, is a linear map from $\mathbb{R}^4$ to $\mathbb{R}$, i.e., $lu = \delta_{t-} u + \frac{\sigma^2(t)}{2} s^2 \delta^2_{s} u \in \text{Hom}(\mathbb{R}^4, \mathbb{R})$.  Also, we can consider the sum of squares in the Tikhonov functional as the 2-norm of a vector: $||{\cdot}||_{2}: \mathbb{R}^{MM} \rightarrow \mathbb{R}$. The facts above imply that we can find a matrix $L \in \text{Hom}(\mathbb{R}^{MM}, \mathbb{R}^{MM})$ that represents the T-shape differential operator; and a corresponding function $J: \mathbb{R}^{MM} \rightarrow \mathbb{R}$ that is ``nearly equivalent'' to $\hat{J}_\beta(u): (\mathbb{R}^M \times \mathbb{R}^{M}) \rightarrow \mathbb{R}$, besides that the domain is different (but feasibly interchangable). This would lead us to the idea of serialization of $u$:  

\begin{equation}
    u_{(jM + i)} = u_{(i, j)}.
\end{equation}

By doing this, we could form a vector $[u_0, u_1, u_2, ..., u_{MM}] \in \mathbb{R}^{MM}$ that helps convert the Tikhonov functional to a much simpler form:
\begin{eqnarray}
    \hat{J}_\beta(u) &=&  dtds \sum_{i \in \hat{Q}_{2\tau}\setminus\partial \hat{Q}_{2\tau}} (lu_{i})^2 + \beta dtds \sum_{i \in \hat{Q}_{2\tau}} (u_i - F_i)^2\\
    &=& dtds ||{Lu}||^{2}_{2} + \beta dtds ||({u - \mathbf{F}})||^{2}_{2}\\
    &=& dtds \left( ||{Lu}||^{2}_{2} + \beta ||({u - \mathbf{F}})||^{2}_{2} \right).
\end{eqnarray}

Where $\mathbf{F} \in \mathbb{R}^{M^2} $ is the vector of values from function $F$ at positions corresponding to $u$. And the next step is to find the matrix $L$. Thanks to the enlightenment from \cite{Wolf}, we found that $L$ can be constructed by Toeplitz-like matrices, diagonal matrices, and the Kronecker product.

Let
\[
  D_t^* = \frac{1}{dt}
  \begin{bmatrix}
  0 & 0 & 0 & 0 & ... & 0\\
  -1 & 1 & 0 & 0 & ... & 0\\
  0 & -1 & 1 & 0 & ... & \vdots\\
  0 & 0 & -1 & 1 & ... & \vdots\\
  \vdots & \ddots & \ddots & \ddots & \ddots & 0\\
  0 & ... & ... & ... & -1 & 1
  \end{bmatrix} \in \mathbb{R}^{M \times M}
\]
and
\[
  D_{ss}^* = \frac{1}{ds^2}
  \begin{bmatrix}
  0 & 0 & 0 & 0 & 0 & ... & 0\\
  1 & -2 & 1 & 0 & 0 & ... & 0\\
  0 & 1 & -2 & 1 & 0 & ... & \vdots\\
  0 & 0 & 1 & -2 & 1 & ... & \vdots\\
  \vdots & \ddots & \ddots & \ddots & \ddots & \ddots & 0\\
  0 & ... & ... & ... & 1 & -2 & 1\\
  0 & ... & 0 & 0 & 0 & 0 & 0
  \end{bmatrix} \in \mathbb{R}^{M \times M}
\]

And the the matrix $L$ can be constructed by 
\[
  L = D_t^* \otimes I^M + R (I^M \otimes D_{ss}^*)
\]
with $R \in \mathbb{R}^{MM \times MM}$ being a diagonal matrix with elements corresponding to the factor $\frac{\sigma^2(t)}{2} s^2$.

\vspace{1em}
\section{Pre-processing boundary and initial conditions}

The boundary and initial conditions needs to be ruled out from the equation $L u = 0$, since these variables are prescribed and cannot be considered as free variables in the system of equations.

\vspace{1em}
\subsection{Why we must remove the boundaries}
\vspace{0.2em}

The necessity of pre-processing the boundary and initial conditions is originated from the nature of differential equations. Consider the original PDE
\[
    \mathcal{L} u \equiv \frac{\partial u}{\partial t} + \frac{\sigma^2(t)}{2} s^2 \frac{\partial^2 u}{\partial s^2} = 0
\]

Suppose $u_0 \in H^{2, 1}(Q_{2\tau})$ satisfies the equation above, then we can easily construct $u_1 = u_0 + C$ as another answer to the PDE, where $C \in \mathbb{R}$ is an arbitrary constant. This phenomenon reflects the truth that the boundary and initial conditions are extremely important to ensure the uniqueness of the final solution.

\subsection{Why Tikhonov-like regularization is not enough}

We have the Tikhonov-like functional
\[
    J_\beta(u) \equiv \int_{Q_{2\tau}} (\mathcal{L}u)^2 dsdt + \beta ||({u - F})||_{L^2(Q_{2\tau})}
\]
to prevent $u$ from deviating too much from $F(s, t)$. Although $F(s, t)$ satisfies the boundary and initial conditions when $(s, t) \in \partial \qtt$, it does not stick the solution to boundary/initial conditions since the regularization parameter $\beta$ is adjustable in running cases. The parameter is fine-tuned during the beta finding session to make sure that we remove the noise from input data as much as possible; then, if the $\beta$ we find eventually is enoughly small, the effectiveness of regularization for constraining the boundary and initial conditions is negligible.

\subsection{Assigning values to free variables in system}

Let's take a small $3 \times 3$ real matrix as a quick example. Consider the system $Ax = b$:

\[
  \begin{bmatrix}
      a & b & c\\
      d & e & f\\
      g & h & i\\
  \end{bmatrix} 
  \begin{bmatrix}
      x_1\\x_2\\x_3
  \end{bmatrix} = 
  \begin{bmatrix}
      b_1\\b_2\\b_3
  \end{bmatrix}
\]

with the following requirements: 

\begin{itemize}
    \item $\dim(\text{null}(A)) = 1$;
    \item $x_1$ is the free variable;
    \item $[a \; d \; g]^T \not= 0$.
\end{itemize}

Then $x_2$ and $x_3$ are linearly dependent on $x_1$, and there exists a unique solution $x = (x_1, x_2, x_3)$ once $x_1$ is given.

Now we define $x_1 = c_1$, and solve the system. 

Firstly, we do finitely many steps of row operations to eliminate any non-zero elements in the first row $[a\;b\;c]$; then we have the new system $A'x = b'$:

\[
  \begin{bmatrix}
      0 & 0 & 0\\
      d' & e' & f'\\
      g' & h' & i'\\
  \end{bmatrix} 
  \begin{bmatrix}
      x_1\\x_2\\x_3
  \end{bmatrix} = 
  \begin{bmatrix}
      0\\b_2'\\b_3'
  \end{bmatrix}
\]
In this case $b'_1$ must equal to $0$, otherwise there will be no solution to the system.

Secondly, we construct an auxillary vector $x_{bd} = (c_1, 0, 0)$ ($bd$ stands for boundary, as the same notation in the PDE case) to eliminate $x_1$ from the system. We substract $A'x_{bd}$ from the RHS to yield the new system $A'x = b' - A'x_{bd}$:
\[
  \begin{bmatrix}
      0 & 0 & 0\\
      d' & e' & f'\\
      g' & h' & i'\\
  \end{bmatrix} 
  \begin{bmatrix}
      x_1'\\x_2\\x_3
  \end{bmatrix} = 
  \begin{bmatrix}
      0\\b_2''\\b_3''
  \end{bmatrix}
\]
Then, the solution to the new system exists only if $x_1' = 0$. This would help us to rewrite $d' = g' = 0$ to the system since they will always multiply by $x_1'$ in the matrix:
\[
  \begin{bmatrix}
      0 & 0 & 0\\
      0 & e' & f'\\
      0 & h' & i'\\
  \end{bmatrix} 
  \begin{bmatrix}
      0\\x_2\\x_3
  \end{bmatrix} = 
  \begin{bmatrix}
      0\\b_2''\\b_3''
  \end{bmatrix}
\]
Finally, shrink the matrix to get our reduced-order system $\overline{A}\overline{x} = \overline{b}$:
\[
  \begin{bmatrix}
      e' & f'\\
      h' & i'\\
  \end{bmatrix} 
  \begin{bmatrix}
      x_2\\x_3
  \end{bmatrix} = 
  \begin{bmatrix}
      b_2''\\b_3''
  \end{bmatrix}.
\]
Suppose the solution to the reduced matrix is $\overline{x_0} = (c_2, c_3)$, then the final answer for the original system, given that $x_1 = c_1$, is $x_0 = (c_1, c_2, c_3)$.

\vspace{1em}

\subsection{Removing boundaries  from $\min(\hat{J_\beta})$ }
\vspace{0.2em}

The finite elements approximation of the Tikhonov-like functional $J_\beta$ is given by
\begin{equation}
    \hat{J_\beta} = ||{Lu}||^{2}_{2} + \beta ||({u - F})||^{2}_{2}
\end{equation}

from above; where $L = D_t + RD_{ss} \in \text{Hom}(\mathbb{R}^{m^2}, \mathbb{R}^{m^2})$. As we can observe from the previous chapter, the rows in matrix $L$ corresponding to the boundaries $\partial Q_{2\tau}$ are all zeros. This means that the matrix $L$ as a nullity of at least $2M - 2$ (the number of points on boundaries). By using the same method mentioned above, we first construct the boundary vector $u_{bd}$, containing given values on the boundaries and zero otherwise; and then subtract the matrix product $L u_{bd}$ from the equation $Lu = 0$ to get
\begin{equation}
  Lu = -L u_{bd}.
\end{equation}

Then we could delete the all-zero rows and columns corresponding to the boundaries to reduce the order of matrix $L$. Notice that $u_{bd} - F_{bd} \equiv 0$ since the $F$ always satisfies the boundary values given; so we can do the same reduction to vector $F$. After the deletion, $L$ becomes  $(M^2 - 2M + 2) \times (M^2 - 2M + 2)$ matrix instead of $M^2 \times M^2$, and $u, F \in \mathbb{R}^{M^2 - 2M + 2}$. We finally produce the reduced minimization problem without changing the solution:
\begin{equation}
    \min \hat{J_\beta} = \min_{u \in \mathbb{R}^{M^2 - 2M + 2}} {\left( ||({Lu})||^{2}_{2} + \beta ||({u - F})||^{2}_{2} \right)}.
\end{equation}

This method ensures that the minimizer $u$ follows strictly with $F$ on the boundaries in order to produce a creditable solution, and slightly reduces the require performance of computer (although the runtime complexity is not reduced).

\vspace{1em}
\section{Minimized $J_\beta$ using conjugate gradient method}

Our ultimate goal is to minimize the function 
\[
  \hat{J}_\beta(u) = {\left( ||({Lu})||^{2}_{2} + \beta ||({u - F})||^{2}_{2} \right)}
\]

In this chapter, Python is used to implement all data structures needed and to minimize the given numerical problem.

\subsection{Normalizing the matrix before conjugate gradient method}

We noticed through practical using of the code, that the matrix $L$ would typically have huge float numbers. It is occurred mainly because the differece between the ask and bid price of the underlying stock is often $0.01$ (or even less in more accurate stock markets). This would result in elements values more than $10^8$ in the matrix $L$, where continuous multiplications ($M^2$ times) of these elements will definitely cause the floating number to overflow (In Python 3, the max precision for float is about $10^{308}$). Therefore, we normalized each row of matrix $L$, in order to prevent this case from happening.

For each row $r_i \in \mathbb{R}^{M^2 - 2M + 2}$ in the system $Lu = b$, we do
\begin{verbatim}
    norm_i = L2norm(r_i)
    r_i = r_i / norm_i
    b_i = b_i / norm_i
\end{verbatim}

This will eliminate the possibility of overflowing while maintaining the correctness of the solution.

\subsection{Minimized using Python SciPy module}

In order to minimize $J_\beta(u)$, conjugate gradient method (CG) is used due to the resemblance of the system of equation $Lu = -Lu_{bd}$ to the symmetric and positive-definite matrix $A$ required by the direct CG: $Ax = b$. However, we chose to use the iterative method of CG since the lost function given by direct CG $f(x) = \frac{1}{2} x^T A x - x^T b$ does not match our given Tikhonov-like functional, which is defined as $J_\beta(u)$.

The SciPy module in Python offers the functionality of minimizing a function using CG. Let \texttt{j\_beta} be a lambda function in Python which takes in a vector of dimension $M^2 - 3M + 2$ and returns a floating number, and \texttt{u} be a mutable array of length $M^2 - 3M + 2$, then calling

\texttt{scipy.minimize(j\_beta, u, method=`CG')}

would return a \texttt{OptimizeResult} containing the minimizer, iterated time, and result status (whether CG succeeds). Extracting the mutated \texttt{u} and reshape it to the matrix of dimension $(M - 2) \times (M - 1)$ would give us the desired output representing the meshgrid of domain $\hat{Q}_{2\tau}$. The results obtained using this method will be shown in the next chapter.

\vspace{1em}

\section{Results}

Tables 1 and 2 display our old results of 368 options we tested in 2015 \cite{KlibGol}.

\vspace{0.5cm
}

\textbf{Table 1. Profits and losses by three different methods} \vspace{0.5cm
}

\begin{tabular}{|l|c|l|}
\hline
Method & Number of tested options & Total profit/loss  \\ \hline
Black-Scholes & $368$ & $+\$11,797,876$ \\ \hline
Last price extrapolation & $368$ & $+\$733,501$ \\ \hline
Ask price extrapolation & $368$ & $-\$505,175,308$ \\ \hline
\end{tabular}

\bigskip
\textbf{Table 2. Percentages of options with profits/losses} \vspace{0.5cm}

\begin{tabular}{|l|l|l|l|}
\hline
Method & profit & loss & zero \\ 
\hline
Black-Scholes & 72.83\% & 16.85\% & 10.32\% \\ \hline
Last price extrapolation & 39.40\% & 56.80\% & 3.80\% \\ \hline
Ask price extrapolation & 10.86\% & 88.34\% & 0.8\% \\ \hline
\end{tabular}

\newpage
Figure 2(a) displays the histogram of profits/loses with the trading strategy devised by Professor Klibanov. Figure 2(b) is a zoomed part of figure 2(a). The profit was $\$13,841,482$   and our loss was $ \$2,043,606. $ The total profit was $\$11,797,876$ (the difference between them).

\begin{figure}[tph]
\centering
\begin{tabular}{cc}
{\includegraphics[width = 0.5\textwidth, height =
0.4\textwidth]{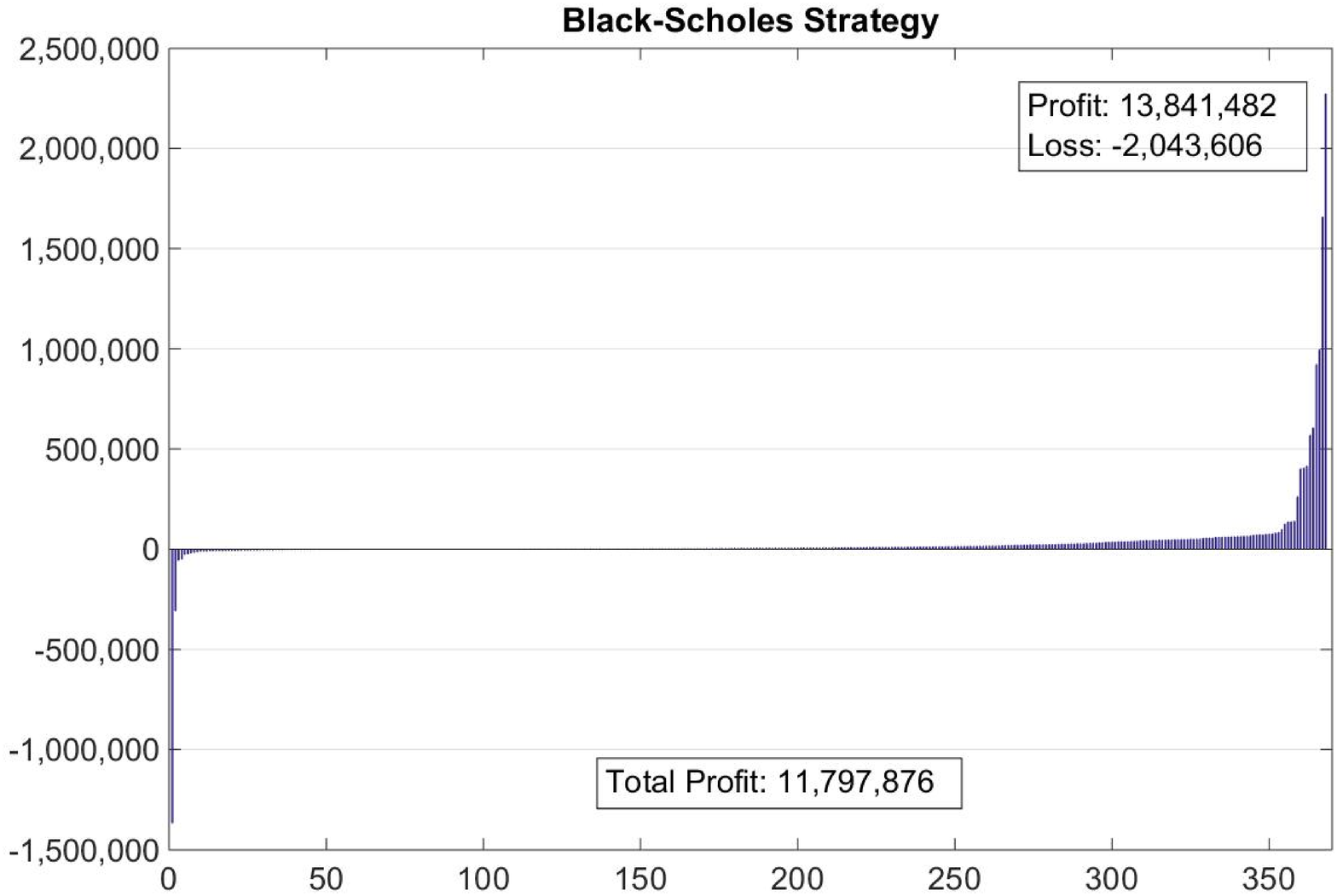}} & {%
\includegraphics[width = 0.5\textwidth,
height = 0.4\textwidth]{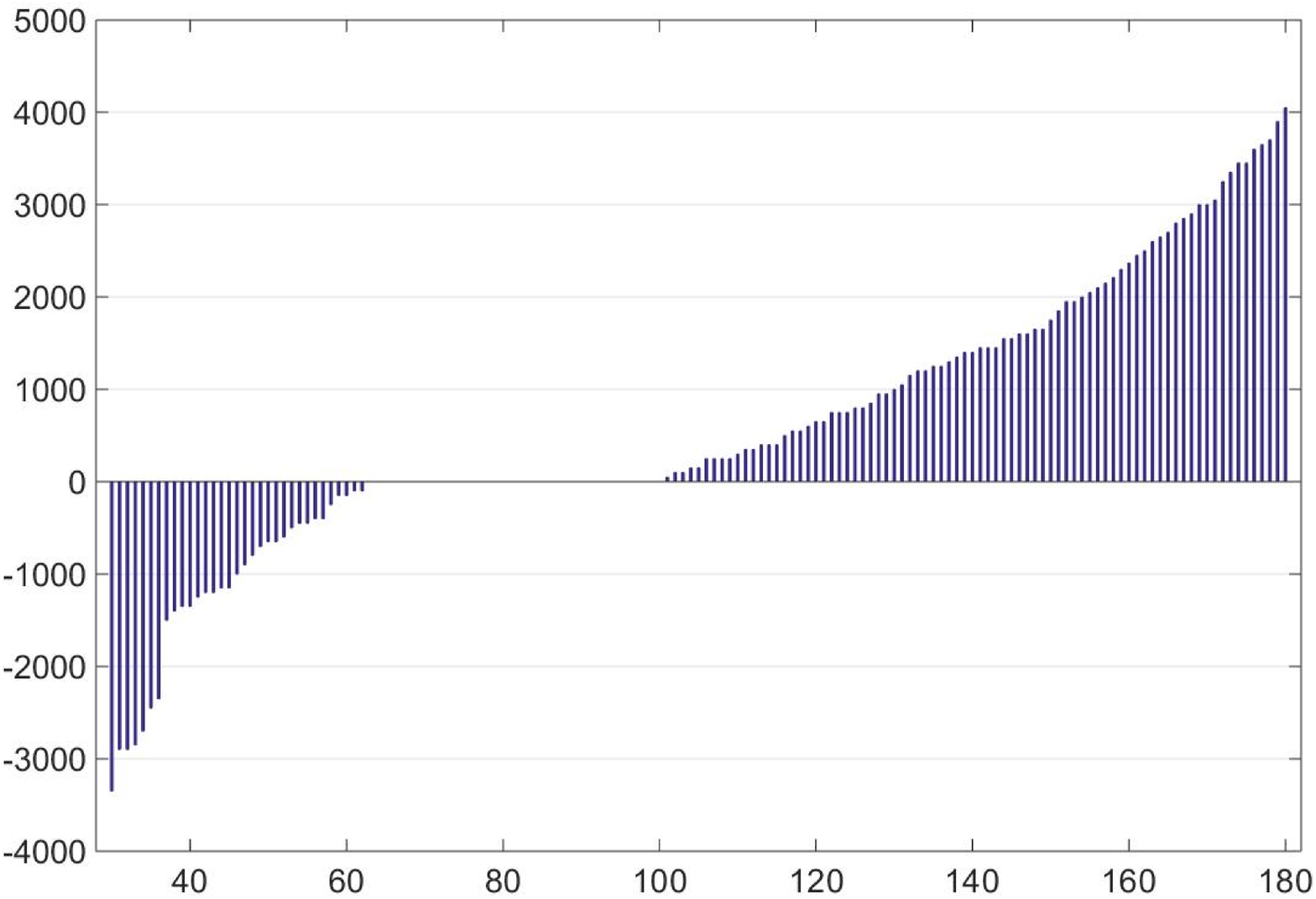}} \\ 
a) & b)%
\end{tabular}
\caption{ The histograms of profits and loses. a) The histogram. b) The zoomed part.
 }
\end{figure}

Figure 3 reflects our new improved results. From the data offered by Bloomberg Terminal, we obtained 50,446 ``data blocks", which consists of option and stock price data for three consecutive trading days (thus Thursday and Friday connected with Monday are considered consecutive). Then we produced the same number of minimizers regarding to the prediction of each of these data blocks. let REAL$_{+1}$ denote the real EOD option price on tomorrow and REAL$_{+2}$ the day after tomorrow, meanwhile EST$_{+1}$ being the estimated option price for tomorrow and EST$_{+2}$ for the day after tomorrow, we define the absolute error:

\[
    err = \frac{1}{2} \left( \frac{\left| \text{EST}_{+1} - \text{REAL}_{+1} \right|}{\text{REAL}_{+1}} + \frac{\left| \text{EST}_{+2} - \text{REAL}_{+2} \right|}{\text{REAL}_{+2}} \right)
\]

In particular, Figure 3(a) displays the histogram of absolute errors of these 50k estimates  compared to real data from Bloomberg terminal. Figure 3(b) represents a zoomed part of figure 3(a).  

\vspace{20em}

\begin{figure}[tph]
\centering
\begin{tabular}{cc}
{\includegraphics[width = 0.5\textwidth, height =
0.4\textwidth]{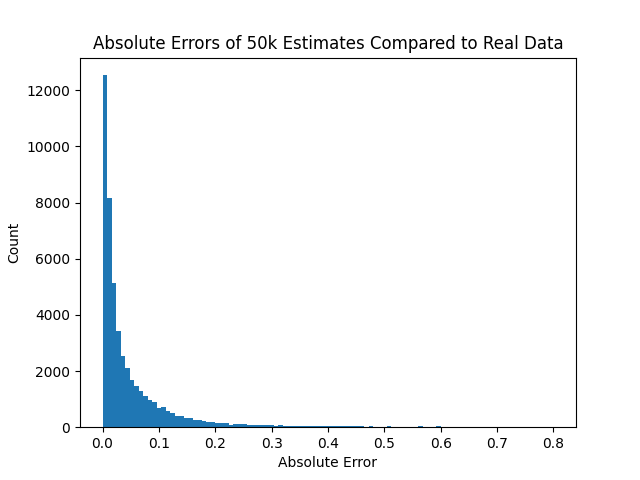}} & {%
\includegraphics[width = 0.5\textwidth,
height = 0.4\textwidth]{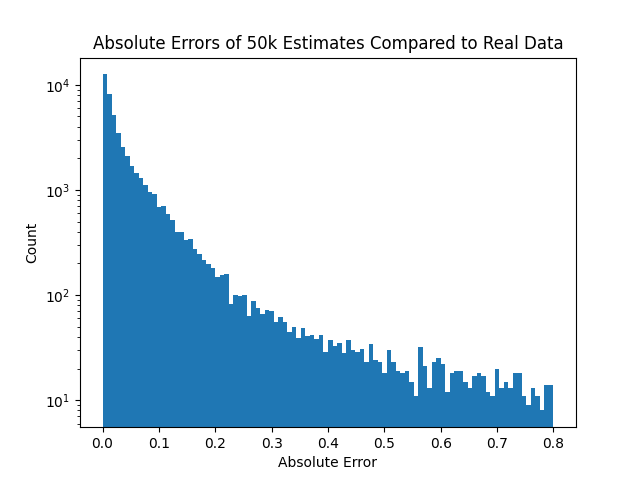}} \\ 
a) & b)%
\end{tabular}
\caption{ The histograms of profits and loses. a) The histogram. b) The zoomed part with logarithm y-axis.
 }
\end{figure}

The same histograms with the best resolution are shown below.
\begin{center}
\includegraphics[scale=0.7]{histogram_50k_log.png}
\includegraphics[scale=0.7]{histogram_50k.png}
\end{center}

From the histogram above, we calculated the median of errors is $2.29\%$.

\section{Application of machine learning to improve profitability of options trading strategy}

\bigskip

\textbf{Method:} We have improved these results using \textit{machine learning} to classify the minimizer's set to filter for the trading strategy inclusion.
We built 13 element input vector consisting of minimizers (for $t=\tau,2\tau$), stock ask and bid  price  (for $t=0$), option ask and bid price and volatility (for $t=-2\tau,-\tau, 0 .$) 
\newpage

Supervised Machine learning has been applied to the binary classification neural network for the logistic probability loss function with regularization: 
\begin{equation}
    L(\theta)=\frac{1}{m}\sum_{i=1}^{m}[-y^{(i)}\log(h_{\theta}(x^{(i)}) -(1-y^{(i)})\log (1-h_{\theta}(x^{(i)}))] + \frac{\lambda}{2m} \sum_{j=1}^{n}\theta_{j}^{2}
\end{equation}

Where $\theta$ are  weights which are optimized by minimizing the loss function using the method of gradient descent.  $\lambda$ is a parameter of regularization. $x^{(i)}$ is our normalized 13 - dimensional vectors. $h_{\theta}$ is output of the neural network. $m$ is the number of vectors in the training set. $y^{(i)}$ 
is our labels (the ground truth). The labels are set to 1 for profitable trades and 0 otherwise.  All vectors and labels are split into three parts: training, verification/validation and test sets. The training set is used for weight learning. Verification/validation set is used for tuning of the neural network hyper-parameters. Test set is for generating the outcomes of trading strategy. We compared the profitability of the trading strategy based on the original minimizer's set with the profitability of the classified set.

\subsection{Results:}

We created a heuristic model for the neural network of three hidden layers of dimension, of 50, 25, 14 dimensions respectively. We also found that learning rate $= 0.00005$ and iterate time $=200$ is a proper parameter settings for this neural network.

Figure 4 shows the learning curve of the neural network, using the MSELoss function with mean reduction as the evaluation loss function.
\begin{figure}[tph]
\begin{center}
    \includegraphics[scale=0.8]{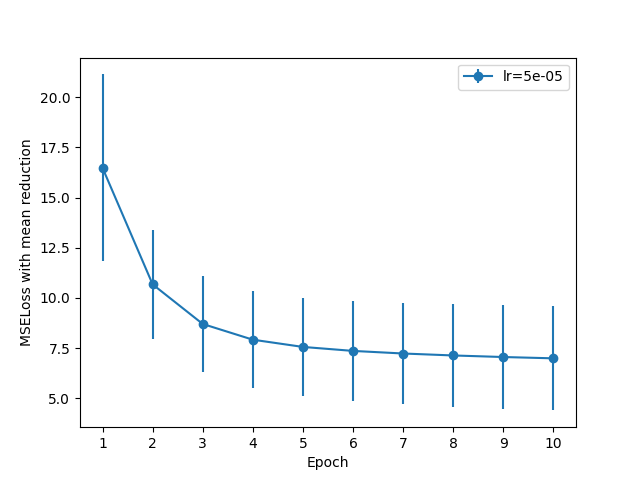}
    \caption{Learning curve from heuristic neural network model}
\end{center}
\end{figure}
We used k-cross validation to ensure the consistency of the neural network, with $k = 10$ in the graph above. The dots in the graph are the mean values of the loss at each epoch (1 epoch $=$ 20 iterations), while the vertical bars represents the standard deviation among the different training datasets. From the graph, we can conclude that the loss converges eventually, and thus machine learning technique is a feasible way to reduce noise and errors produced by pure mathematical solutions.

\newpage
\bibliographystyle{plain}
\bibliography{references}
\end{document}